# COMMENTS ON "A NOTE ON OPTIMAL DETECTION OF A CHANGE IN DISTRIBUTION," BY BENJAMIN YAKIR[1]

BY YAJUN MEI

*Georgia Institute of Technology*

The purpose of this note is to show that in a widely cited paper by Yakir [*Ann. Statist.* **25** (1997) 2117–2126], the proof that the so-called modified Shiryayev–Roberts procedure is exactly optimal is incorrect. We also clarify the issues involved by both mathematical arguments and a simulation study. The correctness of the theorem remains in doubt.

**1. Introduction.** In the change-point literature, as in sequential analysis more generally, theorems establishing exact optimality of statistical procedures are quite rare. Moustakides [2] and Ritov [5] showed that for the simplest problem where both the pre-change distribution $f_0$ and the post-change distribution $f_1$ are fully specified, Page's cumulative sum (CUSUM) procedure [3] is exactly optimal in the sense of minimizing the so-called "worst case" detection delay subject to a specified frequency of false alarms. Earlier, Lorden [1] showed this optimality property holds asymptotically. Besides Page's CUSUM procedure and its generalizations, the most commonly used and studied approach to define change-point procedures is that of Shiryayev [7] and Roberts [6]. Yakir [8] published a proof that claims when both $f_0$ and $f_1$ are fully specified, a modification of the Shiryayev–Roberts procedure is exactly optimal with respect to a slightly different measure of quickness of detection. In this note we show that Yakir's proof is wrong. It is still an open problem whether the modified Shiryayev–Roberts procedure is in fact optimal, although its asymptotic optimality was proved in [4].

**2. Notation.** In this note we use the notation of Yakir [8]. However, there is one ambiguity between $\mathbf{P}_k(\cdot)$ and $\mathbf{P}(\cdot|\nu = k)$ in [8]. The change-point $\nu$ is

Received January 2004; revised July 2005.
[1]Supported in part by the National Institutes of Health under Grant R01-AI055343.

an unknown constant under the non-Bayesian formulation, but it is a random variable in the auxiliary Bayes problem $B(G, p, c)$. To avoid confusion, in this note we denote by $\nu$ the change-point only in the Bayes problem $B(G, p, c)$, and for $1 \leq k \leq \infty$, we denote by $\mathbf{P}_k$ the probability measure (with change time $k$) when the observations $X_1, X_2, \ldots$ are independent such that $X_1, \ldots, X_{k-1}$ have density $f_0$ and $X_k, X_{k+1}, \ldots$ have density $f_1$. In other words, we use $\mathbf{P}_k(\cdot)$ in the context of the non-Bayesian formulation, while $\mathbf{P}(\cdot|\nu = k)$ is used in the context of the Bayes problem $B(G, p, c)$. A critical mistake was made in the proof presented by Yakir [8] because of the confusion between $\mathbf{P}_k(\cdot)$ and $\mathbf{P}(\cdot|\nu = k)$, especially when $k = 1$.

The modified Shiryayev–Roberts procedure, proposed in [4], is defined by

$$N_A^* = \inf\{n \geq 0 : R_n^* \geq A\}, \tag{1}$$

where

$$R_n^* = (1 + R_{n-1}^*)\frac{f_1(X_n)}{f_0(X_n)},$$

and $R_0^* \in [0, \infty)$ has a distribution chosen by the statistician.

For the right distribution of $R_0^*$, the asymptotic optimality of $N_A^*$ was proved in [4]. Later Yakir [8] claimed that $N_A^*$ is exactly optimal in the sense of minimizing the "average" detection delay

$$\mathcal{D}(N) = \sup_{1 \leq k < \infty} \mathbf{E}_k(N - k + 1 | N \geq k - 1) \tag{2}$$

among all stopping times $N$ satisfying $\mathbf{E}_\infty N \geq \mathbf{E}_\infty N_A^*$. In this note, we explain what is wrong with Yakir's proof.

**3. Theoretical results.** In order to prove optimality properties of $N_A^*$ for the right distribution of $R_0^*$, Pollak [4] and Yakir [8] considered the following extended Bayes problem $B(G, p, c)$. Let $G$ be a distribution over the interval $[0, 1]$. Suppose $0 < p < 1$. Assume that a random variable $\pi_0$ is sampled from the distribution $G$ before taking any observations. Given the observed value of $\pi_0$, suppose the prior distribution of the change-point $\nu$ is given by $\mathbf{P}(\nu = 1) = \pi_0$ and $\mathbf{P}(\nu = n) = (1 - \pi_0)p(1-p)^{n-2}$ for $n \geq 2$. Consider the problem of minimizing the risk

$$\mathcal{R}(N) = \mathbf{P}(N < \nu - 1) + c\mathbf{E}(N - \nu + 1)^+,$$

where $c > 0$ can be thought of as the cost per observation of sampling after a change. It is well known [7] that the Bayes solution of this extended Bayes problem $B(G, p, c)$ is of the form

$$M_{G,p,c} = \inf\{n \geq 0 : R_{q,n}^* \geq A\},$$

where $q = 1 - p$, and

$$R^*_{q,0} = \frac{\pi_0 q}{(1-\pi_0)p} - 1, \qquad R^*_{q,n} = (R^*_{q,n} + 1)\frac{g(X_n)}{f(X_n)}\frac{1}{q} \qquad \text{for } n \geq 1,$$

where $\pi_0$ has a distribution $G$. Yakir [8] showed that for some sequence of $p \to 0$, there exists a sequence of $G = G_p$ and $c = c_p$ such that $c \to c^*$ and $\pi_0/p \to R^*_0 + 1$ in distribution, and so $N^*_A$ defined in (1) is a limit of Bayes solutions $M_{G,p,c}$. Yakir [8] claimed that the Bayes solution $M_{G,p,c}$ satisfies (see lines 11–12 on page 2123)

$$(3) \qquad \lim_{p \to 0} \frac{1 - \mathcal{R}(M_{G,p,c})}{p} = (1 - c^*\mathbf{E}_1 N^*_A)[\mathbf{E}R^*_0 + 1 + \mathbf{E}_\infty N^*_A].$$

The proof of the exact optimality of the modified Shiryayev–Roberts procedure in [8] is based on this equation. However, the next theorem shows that equation (3) does not hold in general.

THEOREM 1.
$$(4) \qquad \lim_{p \to 0} \frac{1 - \mathcal{R}(M_{G,p,c})}{p} = [\mathbf{E}R^*_0 + 1 + \mathbf{E}_\infty N^*_A] \\ - c^*[\mathbf{E}_1(R^*_0 N^*_A) + (\mathbf{E}_1 N^*_A)(1 + \mathbf{E}_\infty N^*_A)].$$

PROOF. For the extended Bayes problem $B(G, p, c)$, any stopping rule $N$ satisfies

$$(5) \qquad \frac{1 - \mathcal{R}(N)}{p} = \frac{\mathbf{P}(N \geq \nu - 1)}{p}[1 - c\mathbf{E}(N - \nu + 1|N \geq \nu - 1)].$$

Yakir [8] correctly showed that

$$(6) \qquad \lim_{p \to 0} \frac{\mathbf{P}(N \geq \nu - 1)}{p} = \mathbf{E}R^*_0 + 1 + \mathbf{E}_\infty N^*_A.$$

Arguing as in Lemma 13 of [4], we have

$$(7) \qquad \lim_{p \to 0} \mathbf{E}(M_{G,p,c} - \nu + 1|M_{G,p,c} \geq \nu - 1) \\ = \mathbf{E}_1 N^*_A \frac{\mathbf{E}_\infty N^*_A}{\mathbf{E}R^*_0 + 1 + \mathbf{E}_\infty N^*_A} + \frac{\mathbf{E}R^*_0 + 1}{\mathbf{E}R^*_0 + 1 + \mathbf{E}_\infty N^*_A} \lim_{p \to 0} \mathbf{E}(N^*_A|\nu = 1),$$

and the limiting distribution of $R^*_0$ conditional on $\{\nu = 1\}$ has the density

$$\frac{(x+1)\,d\phi_0(x)}{\int (x+1)\,d\phi_0(x)} = \frac{(x+1)\,d\phi_0(x)}{\mathbf{E}R^*_0 + 1},$$

where $\phi_0(x)$ is the unconditional distribution of $R^*_0$. Yakir [8] made a critical mistake by thinking that the limiting distribution of $R^*_0$ conditional on $\{\nu =$ 1} is just $\phi_0(x)$. Since $R_0^* \geq 0$, the stopping times $N_A^*$ are dominated by the Shiryayev–Roberts stopping time with $R_0^* = 0$. Thus, by the dominated convergence theorem,

$$
\begin{aligned}
\lim_{p \to 0} \mathbf{E}(N_A^* | \nu = 1) &= \lim_{p \to 0} \mathbf{E}(\mathbf{E}(N_A^* | R_0^*, \nu = 1) | \nu = 1) \\
&= \frac{\mathbf{E}_1(N_A^*(R_0^* + 1))}{\mathbf{E} R_0^* + 1}.
\end{aligned}
\tag{8}
$$

Theorem 1 follows at once from (5)–(8) and the fact that $c = c(p) \to c^*$ as $p \to 0$. □

A comparison of equations (3) and (4) shows that the major problem in Yakir's proof comes from the fact that the term $\mathbf{E}_1(R_0^* N_A^*)$ is missing. To further demonstrate this, as suggested by one referee, let us consider

$$
\mathcal{C}(N) = \lim_{p \to 0} \frac{1 - \mathcal{R}(N)}{p}
\tag{9}
$$

for a given stopping time $N$. Since $M_{G,p,c}$ are Bayesian solutions, we have

$$
\lim_{p \to 0} \frac{1 - \mathcal{R}(M_{G,p,c})}{p} \geq \mathcal{C}(N)
\tag{10}
$$

for any given stopping time $N$. Yakir [8] used inequality (10) and equation (3) to prove the exact optimality of $N_A^*$. In the following, we illustrate why Yakir's proof fails.

Note that

$$
\begin{aligned}
1 - \mathcal{R}(N) &= \mathbf{P}(N \geq \nu - 1) - c \mathbf{E}(N - \nu + 1)^+ \\
&= \mathbf{E}_{\pi_0}\left[ \sum_{k=1}^{\infty} P(\nu = k)(\mathbf{P}(N \geq k - 1 | \nu = k, \pi_0) \right. \\
&\quad \left. - c\mathbf{E}((N - k + 1)^+ | \nu = k, \pi_0)) \right],
\end{aligned}
$$

where $\mathbf{E}_{\pi_0}$ denotes expectation with respect to $\pi_0$. Here it is important to point out that $\mathbf{P}(\cdot | \nu = k, \pi_0)$ is same as $\mathbf{P}_k(\cdot)$ but $\mathbf{P}(\cdot | \nu = k)$ is different from $\mathbf{P}_k(\cdot)$ because the prior distribution of $\nu$ depends on $\pi_0$. Since $\mathbf{P}_k(N \geq k - 1) = \mathbf{P}_\infty(N \geq k - 1)$, we have

$$
\begin{aligned}
1 - \mathcal{R}(N) &= \mathbf{E}_{\pi_0}\left[ \sum_{k=1}^{\infty} P(\nu = k)(\mathbf{P}_\infty(N \geq k - 1) - c\mathbf{E}_k(N - k + 1)^+) \right] \\
&= \mathbf{E}_{\pi_0}\left[ \pi_0(\mathbf{P}_\infty(N \geq 0) - c\mathbf{E}_1 N) \right.
\end{aligned}
$$

$$+ (1 - \pi_0) p \sum_{k=2}^{\infty} (1-p)^{k-2}$$

$$\times (\mathbf{P}_\infty(N \geq k-1) - c\mathbf{E}_k(N-k+1)^+)\bigg].$$

Using the facts that $c \to c^*$ and $\pi_0/p \to R_0^* + 1$ in distribution, for *any* given stopping time $N \geq 0$, we have

$$\mathcal{C}(N) = \mathbf{E}_{R_0^*}\bigg[(R_0^* + 1)(\mathbf{P}_\infty(N \geq 0) - c^*\mathbf{E}_1 N)$$

$$+ \sum_{k=2}^{\infty}(\mathbf{P}_\infty(N \geq k-1) - c^*\mathbf{E}_k(N-k+1)^+)\bigg],$$

where $\mathbf{E}_{R_0^*}$ denotes expectations with respect to $R_0^*$. Observe that

$$\mathbf{E}_{R_0^*}((R_0^* + 1)\mathbf{E}_1 N) = \mathbf{E}(\mathbf{E}(R_0^* \mathbf{E}_1 N | R_0^*)) + \mathbf{E}_1 N$$
$$= \mathbf{E}_1(R_0^* N) + \mathbf{E}_1 N$$

because the properties of $R_0^*$ are the same under any probability measure $\mathbf{P}_k$ since $R_0^*$ is chosen by the statistician *before* taking any observations. It is important to point out that $R_0^*$ and the stopping time $N$ *may* or *may not* be correlated under $\mathbf{P}_1$, depending on whether the stopping rule of $N$ involves $R_0^*$. Then, by the facts that $\mathbf{P}_1(N \geq 0) = 1$ and $\sum_{k=2}^{\infty} \mathbf{P}_\infty(N \geq k-1) = \mathbf{E}_\infty N$, we have

$$\mathcal{C}(N) = \mathbf{E} R_0^* + 1 - c^*\mathbf{E}_1(R_0^* N) + \mathbf{E}_\infty N$$

$$- c^* \sum_{k=1}^{\infty} \mathbf{E}_k(N-k+1 | N \geq k-1)\mathbf{P}_\infty(N \geq k-1).$$

On the one hand, for any stopping time $N \geq 0$, by the definition of $\mathcal{D}(N)$ in (2) and the fact that $\sum_{k=1}^{\infty} \mathbf{P}_\infty(N \geq k-1) = \mathbf{E}_\infty N + 1$,

(11) $\quad \mathcal{C}(N) \geq \mathbf{E} R_0^* + 1 - c^*\mathbf{E}_1(R_0^* N) + \mathbf{E}_\infty N - c^*\mathcal{D}(N)(\mathbf{E}_\infty N + 1).$

On the other hand, $N_A^*$ is a so-called *equalizer rule* in the context of the non-Bayesian formulation, that is, for all $k \geq 1$, $\mathbf{E}_k(N_A^* - k + 1 | N_A^* \geq k-1) = \mathbf{E}_1 N_A^* = \mathcal{D}(N_A^*)$. Hence,

(12)
$$\mathcal{C}(N_A^*) = \mathbf{E} R_0^* + 1 - c^*\mathbf{E}_1(R_0^* N_A^*)$$
$$+ \mathbf{E}_\infty N_A^* - c^*\mathcal{D}(N_A^*)(\mathbf{E}_\infty N_A^* + 1),$$

which is exactly the right-hand side of (4) in Theorem 1. Also see Lemma 13 of [4]. Thus, relation (10) is equivalent to stating that $\mathcal{C}(N_A^*) \geq \mathcal{C}(N)$ for any stopping time $N$.

6                                    Y. MEI

Now let us go back the non-Bayesian problem in which we are interested in minimizing the detection delay $\mathcal{D}(N)$ in (2) among all stopping times $N$ satisfying $\mathbf{E}_\infty N \geq \mathbf{E}_\infty N_A^*$. Assume $\mathbf{E}_\infty N_A^* = B$, and let us consider a stopping time $N$ which satisfies the false alarm constraint with equality, that is, $\mathbf{E}_\infty N = B$ (without loss of generality we can limit ourselves to this case). Then from $\mathcal{C}(N_A^*) \geq \mathcal{C}(N)$ and relations (11) and (12), we have

$$\mathbf{E}_1(R_0^* N_A^*) + \mathcal{D}(N_A^*)(1+B) \leq \mathbf{E}_1(R_0^* N) + \mathcal{D}(N)(1+B),$$

from which we *cannot* conclude that

$$\mathcal{D}(N_A^*) \leq \mathcal{D}(N)$$

due to the two terms, $\mathbf{E}_1(R_0^* N_A^*)$ and $\mathbf{E}_1(R_0^* N)$, and the fact that $\mathbf{E}_1(R_0^* N_A^*) \neq \mathbf{E}_1(R_0^*)\mathbf{E}_1(N_A^*)$ since the stopping rule of $N_A^*$ involves $R_0^*$. In [8], the above inequality follows immediately because the two terms, $\mathbf{E}_1(R_0^* N_A^*)$ and $\mathbf{E}_1(R_0^* N)$, are erroneously missing.

**4. Numerical examples.** It is natural to do simulations to confirm that Yakir's result (3) fails while our result (4) is correct. However, it is difficult to simulate the value of the left-hand side of these two equations. Now based on (3), Yakir [8] also showed that

(13)  $$\mathbf{E}_1 N_A^* = \frac{(\mu_0 + 1)(1 - p_0)}{p_0(\mu_0 + 1) + 1},$$

where

$$p_0 = \mathbf{P}(R_0^* \geq A) \quad \text{and} \quad \mu_0 = \mathbf{E}(R_0^* | R_0^* < A).$$

Yakir is correct in deriving (13) as a consequence of (3). Our result (4) and the arguments in [8] lead instead to

(14)  $$\mathbf{E}_1 N_A^* = (\mu_0 + 1)(1 - p_0) - p_0 \mathbf{E}_1(R_0^* N_A^*).$$

Therefore, in order to confirm the incorrectness of Yakir's proof and the existence of the term $\mathbf{E}_1(R_0^* N_A^*)$, it suffices to show that (13) fails while (14) is correct. To illustrate this, we have performed simulations for the following example, which is considered by Pollak [4] and Yakir [8].

Define $f_0(x) = \exp\{-x\}\mathbb{1}(x > 0)$ and $f_1(x) = 2\exp\{-2x\}\mathbb{1}(x > 0)$, and pick an $A$ such that $0 < A < 2$. As shown in [8], the randomized $R_0^* = (R^* + 1)Z$, where $(R^*, Z)$ is uniformly distributed on the set $[0, A] \times [0, 2]$.

It is straightforward to show that $\mu_0 = A/2$, and

$$p_0 = \mathbf{P}(R_0^* \geq A) = 1 - (\log(A+1))/2.$$

Note that Yakir [8] made a minor mistake here by claiming $p_0 = 1 - (\log A)/2$.

Table 1 compares the theoretical values of $\mathbf{E}_1 N_A^*$ given by (13) and (14) to Monte Carlo estimates. Our theoretical result (14) was based on Monte Carlo estimates of $\mathbf{E}_1(R_0^* N_A^*)$, while Yakir's result (13) was calculated exactly. In the Monte Carlo experiment, the number of repetitions was $10^6$ and each result was recorded as the Monte Carlo estimate $\pm$ standard error.

TABLE 1
*Approximations for* $\mathbf{E}_1 N_A^*$

| $A$ | Monte Carlo | Our result (14) | Yakir's result (13) |
|------|---------------------|---------------------|---------------------|
| 1.5  | $0.5799 \pm 0.0007$ | $0.5806 \pm 0.0003$ | 0.4115 |
| 1.6  | $0.6194 \pm 0.0008$ | $0.6197 \pm 0.0004$ | 0.4433 |
| 1.7  | $0.6589 \pm 0.0008$ | $0.6594 \pm 0.0004$ | 0.4757 |
| 1.8  | $0.6993 \pm 0.0008$ | $0.6998 \pm 0.0004$ | 0.5090 |
| 1.9  | $0.7417 \pm 0.0008$ | $0.7404 \pm 0.0004$ | 0.5430 |
| 1.98 | $0.7739 \pm 0.0009$ | $0.7739 \pm 0.0004$ | 0.5708 |

The results in Table 1 suggest that (14) gives correct values for $\mathbf{E}_1 N_A^*$ and (13) does not. These results support the claim that Yakir's proof of exact optimality of the modified Shiryayev–Roberts procedures is flawed.

**Acknowledgments.** The author would like to thank his advisor, Professor Gary Lorden, for his support, guidance, and encouragement. The author would also like to thank the anonymous referees for detailed and helpful suggestions which have improved the presentation of this note.

School of Industrial
  and Systems Engineering
Georgia Institute of Technology
765 Ferst Drive NW
Atlanta, Georgia 30332-0205
USA
E-mail: ymei@isye.gatech.edu